\documentclass[11pt,twoside]{article}

\usepackage{amsmath}
\usepackage{amsthm}
\usepackage{amsfonts, amssymb}
\usepackage{mathrsfs}
\usepackage[all]{xy}
\usepackage{latexsym, stmaryrd}

\def\titulo#1{\noindent{\bf\LARGE{#1}} \bigskip \thispagestyle{plain}}
\def\autor#1{\noindent{\sc #1}\smallskip}
\def\direccion#1{\noindent #1\bigskip}
\def\email#1{\vspace{1cm}\noindent E-mail address: {\sf #1} \bigskip}

\usepackage{latexsym}
\usepackage[dvips]{graphicx}

\theoremstyle{plain}
\newtheorem{lema}{Lemma}[section]
\newtheorem{prop}[lema]{Proposition}
\newtheorem{teo}[lema]{Theorem}
\newtheorem{coro}[lema]{Corollary}
\theoremstyle{remark}

\newtheorem{obs}[lema]{Remark}

\theoremstyle{definition}
\newtheorem{defi}[lema]{Definition}
\newtheorem{ej}[lema]{Example}

\def\s{\mathfrak{S}}

\def\etft0{\mathnormal{etfT_0}}

\def\S{\mathbb{S}}
\def\C{\mathbb{C}}

\begin{document}

\titulo{2-Dimension from the topological viewpoint}

\autor{Jonathan Ariel Barmak, Elias Gabriel Minian}

\direccion{Departamento  de Matem\'atica.\\
 FCEyN, Universidad de Buenos Aires. \\ Buenos
Aires, Argentina}

\begin{abstract}
\noindent In this paper we study the $2$-dimension of a finite poset from the topological point of view. We use homotopy theory of finite topological spaces 
and the concept of a \it beat point \rm to improve the classical results on $2$-dimension, giving a more complete answer to the problem 
of all possible $2$-dimensions 
of  an $n$-point poset.
\end{abstract}



\section{Introduction}

It is well known that finite topological spaces and finite preorders are intimately related. More explicitly, given a finite set $X$, there exists a $1-1$ correspondence between topologies and preorders on $X$ \cite{Ale}. Moreover, $T_0$-topologies correspond to orders.

One can consider this way finite posets as finite $T_0$-spaces and vice versa. Combinatorial techniques based on finite posets together with topological properties can result in stronger theorems \cite{Bar,May,May2,May3,Mcc,Sto}.\\

A basic result in topology says that a topological space $X$ is $T_0$ if and only if it is a subspace of a product of copies of $\s$, the Sierpinski space. 
Furthermore, if $X$ is in addition finite, it is a subspace of a product of finitely many copies of $\s$. 
Using the correspondence between $T_0$-spaces and posets, this result can be expressed as follows: A finite preorder $X$ is a poset if and only if it is a subpreorder of $2^n$ for some $n$.

From this result it seems pretty natural to define the $2$-dimension of a finite $T_0$-space as the minimum $n\in \mathbb{N}_0$ such that $X$ 
is a subspace of $\s ^n$. This coincides with the classical definition of the $2$-dimension of the poset $X$.\\

In 1963, Nov\'ak \cite{Nov} introduced the notion of the $k$-dimension of a finite poset $X$ (for any integer $k\geq 2$), 
extending the definition of dimension of posets given by
Dushnik and Miller \cite{Dus}. One of the first studies on 2-dimension is the foundational paper \cite{Tro} by Trotter. In that paper he introduces the 
notion of the 
$n$-cube $Q_n$ and defines the $2$-dimension of a finite poset $X$ as the smallest positive integer $n$ such that $X$ can be embedded as a subposet of
 $Q_n$. He also proves in \cite{Tro} and \cite{Tro2} the classical formulas and bounds for the $2$-dimension.

In the last 30 years the theory of the 
$2$-dimension was studied by 
many mathematicians and computer scientists. New results and improvements of known results were obtained, but of course there is still much to investigate \cite{Hab,Tro2,Tro3}.\\

In this paper we will show how the classical bounds for the $2$-dimension of a poset of cardinality $n$ can be obtained from the topological 
point of view. Moreover, we will use homotopy theory of finite spaces to improve the classical results on $2$-dimension, giving
 a more complete answer to the problem of all possible $2$-dimensions 
of  $n$-point posets.

The concept of a \it beat point \rm of a finite $T_0$-space (poset) introduced by Stong \cite{Sto} plays an   essential role in our results. 
Explicitly, we prove below the following proposition.

\begin{prop}
Let $X$ be a finite $T_0$-space (poset) and let $x\in X$ be a beat point. Then $$d(X)-1\le d(X\smallsetminus \{x\})\le d(X)$$
\end{prop}

Here $d(X)$ denotes the $2$-dimension of $X$. This result improves (in the case of beat points) the \it continuity \rm property of the $2$-dimension. As 
an immediate consequence we have:

\begin{coro} 
Let $X$ be a finite contractible $T_0$-space. Then $d(X)\le |X|-1$.
\end{coro}

Using these results and the notion of \it non-Hausdorff suspension \rm of a topological space \cite{Mcc}, we deduce our main theorem:

\begin{teo}
Given $n\ge 2$ and $m$ such that $\lceil log_2 n \rceil \le m\le n$, there exists a $T_0$-space (poset) $X$ of cardinality $n$ with $d(X)=m$. 
Moreover, if $m\neq n$, $X$ can be taken contractible.
\end{teo}

%
%

\section{Preliminaries: topologies, preorders and initial maps}

We start by recalling the basic correspondence between topologies and preorders on a finite set.

The first mathematician who related finite topological spaces with preorders was Alexandroff \cite{Ale}. Many years later, Stong \cite{Sto} 
and McCord \cite{Mcc} continued Alexandroff's ideas. Recently, a paper by Osaki \cite{Osa} and a beautiful series of notes by Peter May \cite{May,May2,May3}
 captured the attention of algebraic topologists. In \cite{Bar} we used combinatorial techniques based on finite spaces and posets to solve topological 
 problems, concerning the homotopy groups of the spheres.
 
\bigskip

Given a finite topological space $(X,\tau)$ and $x\in X$, we define the \textit{minimal open set} $U_x$ of $x$ as the intersection of all the open sets containing $x$. The preorder on $X$ associated to the topology $\tau$ is defined as follows: $x\le y$ if $x\in U_y$.\\

Conversely, given a preorder $\le$ on $X$, we define for each $x\in X$ the sets $$U_x=\{y\in X \ | \ y\le x\}$$
It is easy to see that these sets form a basis for a topology, which will be the topology associated to $\le$.

These applications are mutually inverse. Moreover, the order relations on $X$ are in co\-rres\-pon\-dence with the $T_0$-topologies on $X$. Recall that a space $X$ is said to be $T_0$ if for every pair of points in $X$ there exists some open set that contains one and only one of those points.\\

From now on, we will identify finite $T_0$-spaces with finite posets.

\begin{ej}
Let $X=\{a,b,c,d\}$ the 4-point space whose open subsets are $\emptyset, \{a,c,d\},$ $\{b,d\}, \ \{c,d\}, \ \{d\}, \ \{a,b,c,d\}, \ \{b,c,d\}$. Then $X$ is a $T_0$-space with Hasse diagram

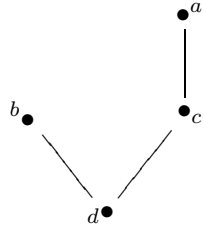
\begin{displaymath}
 \xymatrix@C=10pt{     	&         &      \ \bullet ^a \ar@{-}[d] \\
			^b \bullet \ \ar@{-}[dr] &  & \ \bullet _c \ar@{-}[dl] \\
			&         _d \bullet \   &   }
 \end{displaymath}
 
\end{ej}

Note that there exists a bijection between the topology of a $T_0$-space $X$ and the antichains of the poset $X$ that assigns to each open subset of $X$ the antichain of its maximal elements.\\

\begin{obs} \label{preserva}
The correspondence between topologies and preorders takes products to products, disjoint unions to disjoint unions and subspaces to subpreorders.
\end{obs}

It is not hard to prove that a function between finite spaces is continuous if and only if it is order preserving.\\

The concept of initial map or initial topology is related to the notions of subspace and product \cite{Bou}. 

A map $f:X\to Y$ is initial if the topology on $X$ is induced by $f$. More precisely

\begin{defi}
A function $f:X\to Y$ between topological spaces is an \textit{initial map} (or $X$ has the \textit{initial topology} with respect to $f$) if the topology of $X$ is the coarsest such that $f$ is continuous. Explicitly, $U\subseteq X$ is open if and only if there exists an open set $V\subseteq Y$ such that $U=f^{-1}(V)$.

More generally, we say that a family $\{f_{\lambda}:X\to X_{\lambda}\}_{\lambda \in \Lambda}$ of functions between topological spaces is an \textit{initial family} (or that $X$ has the \textit{initial topology} with respect to $\{f_{\lambda}\}_{\lambda \in \Lambda}$) if the topology of $X$ is the coarsest such that $f_{\lambda}$ is continuous for every $\lambda \in \Lambda$.
\end{defi}

\begin{prop} \label{iniequi}
Let $\{f_{\lambda}:X\to X_{\lambda}\}_{\lambda \in \Lambda}$ be a family of functions between topological spaces. The following are equivalent
\begin{enumerate}
\item[$(i)$] $\{f_{\lambda}\}_{\lambda \in \Lambda}$ is initial.
\item[$(ii)$] $\{f_{\lambda}^{-1}(U) | \ \lambda \in \Lambda , \ U\subseteq X_{\lambda}$ is open$\}$ is a subbase of the topology of $X$.
\item[$(iii)$] For every space $Z$ and every function $g:Z\rightarrow X$, $g$ is continuous if and only if $f_{\lambda}g$ is continuous for every $\lambda \in \Lambda$.
\end{enumerate} 
\end{prop}

\begin{ej}
Let $\{X_{\lambda}\}_{\lambda \in \Lambda}$ be a family of spaces. Then the family of projections $\{ \underset{\lambda \in \Lambda}{\prod} X_{\lambda} \overset{p_{\gamma}}{\to}X_{\gamma} \}_{\gamma \in \Lambda}$ is an initial family.
\end{ej}

It is easy to see that a family $\{f_{\lambda}:X\to X_{\lambda}\}_{\lambda \in \Lambda}$ is initial if and only if the induced function $f:X\to \underset{\lambda \in \Lambda}{\prod} X_{\lambda}$ is an initial map.\\

If $A\subseteq Y$ is a subspace, then the inclusion $i:A\hookrightarrow Y$ is initial. One can generalize the concept of subspace as follows.

\begin{defi}
A function $i:X\to Y$ is a \textit{subspace map} if it is initial and injective.
\end{defi}

If $i:X\to Y$ is a subspace map, $i$ becomes a homeomorphism between $X$ and its image $i(X)$ viewed as a subspace of $Y$. We simply denote $X\subseteq Y$.\\

We denote $\s=\{0,1\}$ the Sierpinski space, whose unique proper open set is $\{0\}$. Note that $\s=2$ is the finite $T_0$-space (order), with $0< 1$.\\

As we pointed out in the introduction, it is a basic topological fact that a space is $T_0$ if and only if it is a subspace of a product of copies of $\s$ (cf. \cite{Bou}). Explicitly

\begin{prop} \label{conocido}
Let $X$ be $T_0$-space. Then the function $$i:X \rightarrow \underset{ \textrm{continuous}}{\underset{h:X\rightarrow \s}{\prod}} \s $$ defined by $i(x)=(h(x))_h$ is a subspace map.
\end{prop}

We will use this result to define the $2$-dimension of a finite $T_0$-space.

\section{Initial maps and subposets}

We give a characterization of inicial maps between finite spaces in terms of preorders.

\begin{prop} \label{inicial}
Let $i:X\to Y$ be a function between finite spaces. The following are equivalent:
\begin{enumerate}
\item[$(i)$] $i$ is initial.
\item[$(ii)$] For every $x,x'\in X$ it holds that $x\le x'$ if and only if $i(x)\le i(x')$.
\end{enumerate}
\end{prop}

\begin{proof}
If $i$ is an initial map, then it is continuous and therefore order preserving. Suppose that $x,x'\in X$ are such that $i(x)\le i(x')$, and suppose that $i^{-1}(U)$ is an open set of $X$ that contains $x'$. Then $i(x')\in U$, and therefore $i(x)\in U_{i(x')}\subseteq U$ which implies that $x\in i^{-1}(U)$. It follows that $x\in U_{x'}$ and hence $x\le x'$.\\

Conversely, if condition $(ii)$ holds, $i$ is order preserving and then continuous. We want to show that the topology $\tau$ of $X$ is the coarsest that makes $i$ continuous. Suppose that $\tau '$ is another topology which induces the preorder $\preceq$ and makes $i$ continuous. If $x\preceq x'$, then $i(x)\le i(x')$ and, by $(ii)$, $x\le x'$. Since $id:(X,\tau')\to (X,\tau)$ is order preserving, it is continuous, and then $\tau \subseteq \tau'$. 
\end{proof}

Suppose now that $i:X\to Y$ is initial and $x,x'$ are such that $i(x)=i(x')$. Since $i(x)\le i(x')$ and $i(x)\ge i(x')$, it follows by the previous proposition that $x\le x'$ and $x\ge x'$. Therefore, if the preorder on $X$ is antisimetric, $i$ results injective. Then we obtain the following

\begin{coro} \label{bastaini}
Let $X, \ Y$ be finite spaces such that $X$ is $T_0$ and suppose that $i:X\rightarrow Y$ is an initial map. Then $i$ is a subspace map.
\end{coro}

\begin{obs}
In fact, the previous result works for infinite spaces as well. But, as we have seen, the proof is very simple in the finite space case by the tractability of posets.
\end{obs}

By \ref{conocido} every finite $T_0$-space is a subspace of a product of finitely many copies of $\s$. Then, the following definition makes sense.

\begin{defi}
Let $X$ be a finite $T_0$-space. We define the \textit{$2$-dimension} of $X$ as the minimum $n\in \mathbb{N}_0$ such that $X$ is a subspace of a product of $n$ copies of $\s$.
\end{defi}

Note that the $2$-dimension is well defined by \ref{conocido} and that $2^{|X|}$ is an upper bound. Here $|X|$ denotes the cardinality of $X$. The $2$-dimension of $X$ will be denoted $d(X)$.

\begin{obs}
By \ref{preserva} it is clear that $X$ is a subspace of a product of $n$ copies of $\s$ if and only if $X$ is a subposet of $2^n$. From \ref{conocido} one can deduce that every finite poset is a subposet of a finite boolean algebra. It is easy to see that our definition of $2$-dimension coincides with the classical definition in terms of posets.
\end{obs}

A topological interpretation of the monotony of the $2$-dimension could be the following. If $X$ is a subposet of $Y$, then $X$ is a subspace of $Y$, which is a subspace of $\s ^{d(Y)}$. Therefore $d(X)\le d(Y)$.\\

Now we prove the well-known result  on the  bounds of the $2$-dimension of an $n$-point $T_0$-space, in terms of topology.

\begin{prop} \label{cotas}
Let $X$ be a finite non-empty $T_0$-space. Then $$\lceil log_2 |X|\rceil \le d(X)\le |X|$$
\end{prop}

\begin{proof}
The first inequality is trivial because $X\subseteq \s ^{d(X)}$ and then $|X| \le 2^{d(X)}$. In order to prove the second inequality let us consider the function $h:X\rightarrow \underset{x\in X}{\prod} \s=\s^{|X|}$ defined by $h(y)=(\chi_{U_x^c}(y))_{x\in X}$, where $\chi_{U_x^c}$ is the characteristic function of $U_x^c$ (i.e. $\chi_{U_x^c}(y)=1$ if $y\in U_x^c$ and $\chi_{U_x^c}(y)=0$ in other case).

This is a continuous map because $p_xh=\chi_{U_x^c}$ is continuous for every $x\in X$.

If $U\subseteq X$ is open, $$U=\underset{x\in X}{\bigcup} U_x= \underset{x\in X}{\bigcup} \chi_{U_x^c}^{-1}(\{0\})= \underset{x\in X}{\bigcup} h^{-1}(p_x^{-1}(\{0\}))=h^{-1} (\underset{x\in X}{\bigcup} p_x^{-1}(\{0\}))$$

Since $p_x^{-1}(\{0\})\subseteq \s^{|X|}$ is open for every $x\in X$, $U$ is an open set in the initial topology. It follows that $h$ is initial, 
and by \ref{bastaini}, $h$ is subspace.
\end{proof}

\section{Homotopy and 2-dimension}

In 1966 Stong classified finite spaces by their homotopy types \cite{Sto}. His ideas turned out to be very illuminating to improve the well known result of \ref{cotas}.\\

First we recall some definitions from \cite{Sto}, although this formulation is from \cite{May}.

\begin{defi} 
Let $X$ be a finite $T_0$-space. We say that $x\in X$ is an \textit{up-beat point} if $\{y\in X | \ y>x\}$ has a minimum. Analogously, we say that $x$ is a \textit{down-beat point} if $\{y\in X | \ y<x\}$ has a maximum. In either of these cases we say that $x$ is a \textit{beat point}.

We say that a finite space is a \textit{minimal finite space} if it is $T_0$ and it has no beat points.
\end{defi}

We recall that a subspace $X$ of a topological space $Y$ is a strong deformation retract of $Y$ if there exists a continuous retraction $r:Y \to X$ of the inclusion $i:X \hookrightarrow Y$, such that $ir$ and the identity of $Y$ are homotopic with a homotopy that is stationary on $i(X)$.

Note that if $X$ is a strong deformation retract of $Y$, $X$ and $Y$ have the same homotopy type.\\

Stong proves that if $x$ is a beat point of a finite $T_0$-space $X$, then $X\smallsetminus \{x\}$ is a strong deformation retract of $X$.

\begin{defi}
Let $X$ be a finite space. A subspace $Y\subseteq X$ is a \textit{core} of $X$ if it is a minimal finite space which is a strong deformation retract of $X$.
\end{defi}

Every finite space $X$ has a core, the core $X_c$ is unique up to isomorphism and it is the smallest space which is homotopy equivalent to $X$. Moreover, if $X$ is $T_0$, there exists a sequence $X=X_0\supset X_1 \supset X_2 \supset \ldots \supset X_n=X_c$ where $X_{i+1}$ is obtained from $X_i$ by removing a beat point.\\

\begin{obs} \label{maximo}
If $X$ is a finite $T_0$-space with maximum $m$, a maximal point $x$ of $X\smallsetminus \{m\}$ is an up-beat point. 
Now $X\smallsetminus \{x\}$ is homotopy equivalent to $X$ and has maximum $m$. By an inductive argument, it follows that $X$ is homotopy equivalent to a point, i.e. contractible. 
A topological proof of this fact can be found in \cite{May,Mcc}.
\end{obs}

\smallskip

If $X$ is a finite space, we can define the space $X^{op}$ whose open sets are the closed sets of $X$. It is easy to prove that the induced preorder in $X^{op}$ is the opposite of $X$.\\

It is well known that if $X$ is a finite poset and $x\in X$, then $d(X)-2\le d(X\smallsetminus \{x\})\le d(X)$. 
We improve this result in the case that $x$ is a beat point.

\begin{prop} \label{continuidad}
Let $X$ be a finite $T_0$-space (poset) and let $x\in X$ be a beat point. Then $$d(X)-1\le d(X\smallsetminus \{x\})\le d(X)$$
\end{prop}

\begin{proof}
The second inequality is clear by the monotony of $d$.

Now suppose that $x\in X$ is an up-beat point. Then there exists $y>x$ such that $z>x$ implies $z\ge y$. Let $i:X\smallsetminus \{x\}\hookrightarrow \s^{d(X\smallsetminus \{x\})}$ be a subspace map. In this case we define $i':X\rightarrow \s^{d(X\smallsetminus \{x\})+1}$ in the following way
\begin{displaymath}
  i'(z)=\left\{\begin{array}{ll}
  i(y)0 & \textrm{if } z=x \\
  i(z)0 & \textrm{if } z<x \\
  i(z)1 & \textrm{in other case} \\
 \end{array} \right.
\end{displaymath}  

Here, we regard the elements of $\s ^{d(X\smallsetminus \{x\})+1}$ as $(d(X\smallsetminus \{x\})+1)$-tuples of $0$'s and $1$'s. 
By $i(y)0$ we mean the $d(X\smallsetminus \{x\})$-tuple $i(y)$ followed by a $0$.

If we prove that $i'$ is an initial map, then by \ref{bastaini}, $i'$ is subspace and $d(X)\le d(X\smallsetminus\{x\})+1$.

\smallskip

We will use proposition \ref{inicial} to show that $i'$ is initial. 
We first show that $i'$ is order preserving. Suppose $z<z'$

\begin{itemize}
	\item If $z=x$, $z'>x$ and then $z'\ge y$. Since $i$ is order preserving, $i(z')\ge i(y)$. Therefore $i(z')1\ge i(y)0$ and $i'(z')\ge i'(z)$.
	\item If $z'=x$, $z<x<y$, then $i(z)\le i(y)$ and therefore $i(z)0\le i(y)0$. It follows that $i'(z)\le i'(z')$.
	\item The case $z\neq x \neq z'$ is clear.

%
\end{itemize}

Now suppose that $z\neq z'$ are such that $i'(z)\le i'(z')$.
\begin{itemize}
	\item If $z=x$, $i(y)0\le i'(z')$ and then $i(y)\le i(z')$. By \ref{inicial} $z=x<y\le z'$.
	\item If $z'=x$, $i'(z')$ ends in $0$ and then, so does $i'(z)$. Therefore $z\le x=z'$.
	\item In other case, $i'(z)\le i'(z')$ implies $i(z)\le i(z')$ and since $i$ is initial, $z\le z'$.
\end{itemize}
Again by \ref{inicial}, it follows that $i'$ is initial.\\

If $x\in X$ is a down-beat point, the result follows from the up-beat point case, considering $X^{op}$.
\end{proof}

\begin{coro} \label{contractil}
Let $X$ be a finite contractible $T_0$-space. Then $d(X)\le |X|-1$.
\end{coro}
\begin{proof}
Since $X$ is contractible, its core $X_c$ consists of only one point. If $X=X_c=\raisebox{1pt}{$*$}$ the result is trivial. Otherwise, $X$ is not minimal and it has a beat point $x$. The space $X\smallsetminus \{x\}$ is a strong deformation retract of $X$ and then contractible. By an inductive argument $$d(X\smallsetminus \{x\})\le |X\smallsetminus \{x\}|-1=|X|-2$$
Now the result follows from the previous proposition.
\end{proof}

There are two standard constructions in Topology, more precisely in Homotopy Theory: the cone and the suspension of a topological space. These 
constructions are not useful for finite spaces, but in the finite case there are two analogous constructions introduced in \cite{Mcc}, namely,  
the \textit{non-Hausdorff cone} and the \textit{non-Hausdorff suspension}.

\begin{defi} (McCord)
We define $\C (X)$ the \textit{non-Hausdorff cone} of a space $X$ as the space $X\cup \{*\}$, whose open sets are the open sets
 of $X$ together with $X\cup \{*\}$.

We define $\S (X)$ the \textit{non-Hausdorff suspension} of $X$ as the space $X\cup \{+,-\}$, whose open sets are those of $X$, together with 
$X\cup \{+\}$, $X\cup \{-\}$ and $X\cup \{+,-\}$. We define recursively the \textit{n-fold non-Hausdorff suspension} of $X$ by $\S ^{n}(X)=\S (\S^{n-1}(X))$.
\end{defi}

Note that if $X$ is $T_0$ and finite, the induced order on $\C (X)$ is the one that results when we add a maximum to $X$, and $\S (X)$  results when we add to $X$ two incomparable points $+$ and $-$ which are greater than all the points of $X$. 

In terms of the join operation $\raisebox{1pt}{$\oplus $}$ of posets, it is easy to see that $\C (X)=X \raisebox{1pt}{$\oplus \; \! * $}$ and $\S (X)=X \raisebox{1pt}{$\oplus $} S^0$, where $\raisebox{1pt}{$*$}$ is the $1$-point poset and $S^0=D_2$ denotes the $2$-point antichain ($2$-point discrete space).

\begin{ej}
Let $X=\s \sqcup \raisebox{1pt}{$*$}$, then the Hasse diagrams of $X$, $\C (X)$ and $\S(X)$ are the following
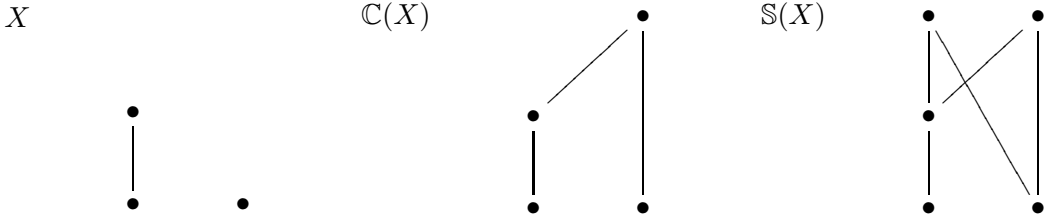
\begin{displaymath}
 \xymatrix@C=30pt{ X & & \\
 		& \bullet \ar@{-}[d] & & \\
		& \bullet & \bullet & }
 \xymatrix@C=30pt{ \C(X) & & \bullet \ar@{-}[dl] \ar@{-}[dd] \\
 		& \bullet \ar@{-}[d] & & \\
		& \bullet & \bullet & 	   }
 \xymatrix@C=30pt{ \S(X) & \bullet \ar@{-}[d] \ar@{-}[ddr] & \bullet \ar@{-}[dl] \ar@{-}[dd] \\
 		& \bullet \ar@{-}[d] & & \\
		& \bullet & \bullet & 	   }
 \end{displaymath}
\end{ej}

The construction of $\S (X)$ becomes important when we work with finite models of spheres (finite spaces with the same homotopy groups of the spheres) \cite{Bar,Mcc}. However we give here a completely different use of this construction.

We deduce from Theorem 2.3 of \cite{Tro} that for every finite $T_0$-space $X$, $$d(\S (X))=d(X)+d(S^0)=d(X)+2$$ 
Inductivelly one proves that $d(\S^n(S^0))=2n+2=|\S^n(S^0)|$ and $d(\S^n(D_3))=2n+3=|\S^n(D_3)|$ for every $n\ge 0$. Therefore we have the following

\begin{lema} \label{cotasup}
Given $n\ge 2$, there exists a $T_0$-space (poset) with $n$ points and whose $2$-dimension is $n$. 
\end{lema}

The next result follows from \ref{maximo} and \ref{contractil}.

\begin{lema} \label{lineal}
Let $n\in \mathbb{N}$ and $X=\llbracket 1,n \rrbracket$ the $T_0$-space with the usual order. Then $d(X)=n-1$. 
\end{lema}

Now we state the main result of this paper.

\begin{teo}
Given $n\ge 2$ and $m$ such that $\lceil log_2 n \rceil \le m\le n$, there exists a $T_0$-space (poset) $X$ of cardinality $n$ with $d(X)=m$. 
Moreover, if $m\neq n$,  $X$ can be taken contractible.
\end{teo}
\begin{proof}
The case $m=n$ follows from \ref{cotasup}, so it remains to analize the case $\lceil log_2 n \rceil \le m\le n-1$.

Let $X_1=\{u_1,\ldots,u_n\}$ be an $n$-point subspace of $\s^{\lceil log_2n \rceil}$ such that $u_1$ is the maximum of $\s^{\lceil log_2n \rceil}$. By \ref{cotas}, it follows that 
\begin{equation} \label{uno}
d(X_1)= \lceil log_2n \rceil
\end{equation}
We define recursively $X_2, \ldots ,X_n$ by $X_{i+1}=\mathbb{C}(X_i\smallsetminus\{u_{i+1}\})$ and we denote by $u_i'$ the maximum of $X_i$. This means that the spaces are constructed from the previous ones by removing a point and putting it in the top.

For every $2\le i \le n$, $u_i'$ is a down-beat point of $X_i$ because $X_i\smallsetminus \{u_i'\}$ has maximum. By \ref{continuidad}, if $1\le j <n$, $d(X_{j+1}\smallsetminus\{u_{j+1}'\})\ge d(X_{j+1})-1$. However, $X_{j+1}\smallsetminus\{u_{j+1}'\}=X_j\smallsetminus\{u_{j+1}\}$ is a subspace of $X_j$, and then $d(X_j)\ge d(X_{j+1}\smallsetminus\{u_{j+1}'\})\ge d(X_{j+1})-1$. Therefore 
\begin{equation} \label{dos}
d(X_{j+1})\le d(X_j) +1
\end{equation}

Every space in this sequence has maximum, and by \ref{maximo} is contractible. Moreover, $X_n$ is linear, therefore, by \ref{lineal} 
\begin{equation} \label{tres}
d(X_n)=n-1
\end{equation}

From (\ref{uno}), (\ref{dos}) and (\ref{tres}) it follows that the sequence $d(X_i)$, with $1\le i\le n$, takes all the values between $\lceil log_2n \rceil$ and $n-1$. This completes the proof.
\end{proof}

\email{jbarmak@dm.uba.ar, gminian@dm.uba.ar}

\end{document}